\documentclass{article}
\usepackage{mathrsfs}
\usepackage[leqno]{amsmath}
\usepackage{amsthm}
\usepackage{amssymb}
\usepackage[all]{xy}
\usepackage[colorlinks,linkcolor=black,citecolor=green,urlcolor=red]{hyperref}
\usepackage[margin=10pt,font=small,textfont=sl,labelfont=bf]{caption}
\newtheoremstyle{mythm}{2ex plus 1ex minus .2ex}{2ex plus 1ex minus .2ex}     {\itshape}{}{\bfseries}{.}{0.7em}{}
\theoremstyle{mythm}
\newtheorem{prop}{Proposition}[section]
\newtheorem{thm}[prop]{Theorem}
\newtheorem{lem}[prop]{Lemma}
\newtheorem{claimm}[prop]{claim}

\newtheoremstyle{mythm1}{2ex plus 1ex minus .2ex}{2ex plus 1ex minus .2ex} {\normalfont}{}{\bfseries}{.}{1em}{} 
\theoremstyle{mythm1}

\newtheorem*{rk}{Remark}
\newtheoremstyle{mythm2}{1.5ex plus 1ex minus .2ex}{1.5ex plus 1ex minus .2ex}     {\normalfont}{\parindent}{\bfseries}{}{1em}{}
\theoremstyle{mythm2}
\newtheorem*{abs}{\textbf{Abstract}}
\usepackage[b5paper,text={125mm,195mm},centering]{geometry}
\numberwithin{equation}{section}
\begin{document}
\title{\textbf{Central Limit Theorems in Deterministic Systems}}
\author{Yuwen Wang, Fudan University}
\date{}
\maketitle
\begin{abs}
    This is a note on some results of the central 
    limit theorem for deterministic dynamical 
    systems. First, we give the central 
    limit theorem for martingales, which is a main tool. Then we give the main results on 
    the central limit theorem in dynamic system in the cases of
    martingale and backward martingale.\\

\end{abs}

\section{Introduction}
I learned the central limit theorem in dynamic systems under 
the guidance of Teacher Xie Jiansheng. This is a note taking down 
the main theorems and conclusions in the literature 
I studied. These results come from  \cite {1},  \cite {2},  
\cite {3},  \cite {7} in the references. In order to make 
it easy for readers to understand, I may revise the original 
literature proof, and I will also use the results in books  
\cite {4},  \cite {5},  \cite {6} to support the 
corresponding proof.

Without specification, the proof of the results in
 this paper is always based on the following assumption:
 there is a probability space
$(\Omega, \mathcal{F}, P)$ and a measure preserving mapping of
 $T:\Omega\rightarrow\Omega$, i.e. for $\displaystyle
A\in \mathcal{F}$, then $P(A)=P(T^{-1}(A))$. It is also 
assumed that if $A\in\mathcal{F}$, then 
$T(A)\in \mathcal{F}$, and the measure preserving system
 $(P,T)$ are ergodic. Suppose $X$ is a measurable 
 function on $\Omega$, denote $U:X\mapsto U(X)=X\circ T$.

 Under these assumptions, for the measurable function $f$ 
 on $\Omega$ which satisfies some conditions, 
 we expect the following to be true: 
\begin{equation*}
\frac{1}{\sqrt{n}}\sum_ {i=0}^{n-1} U^if\rightarrow N (0, \sigma^2)
\end{equation*}
To attain this goal, we will use method of martingale approximation. 
Set $\mathcal{F}_ 0\subset \mathcal{F}$ is a sub $\sigma$ 
algebra, and $\mathcal{F}_ k=T^{-k}\mathcal{F}_ 0, k\in \mathbb{Z}$, 
then there are two general cases at this time
\begin{equation}\label{jiashe}
\begin{aligned}
(1) &\dots\supset\mathcal{F}_ {-1}\supset\mathcal{F}_ 0\supset\mathcal{F}_ 1\supset\dots\\
(2) &\dots\subset\mathcal{F}_ {-1}\subset\mathcal{F}_ 0\subset\mathcal{F}_ 1\subset\dots
\end{aligned}
\end{equation}
In fact, (1) and (2) correspond respectively to 
inverted martingales and martingales. 
Section 2 will describe the central limit theorem of 
martingales. Results of inverted martingales and martingale
are given respectively in section 3 and section 4.

\begin{rk}\label{zhu1}
\begin{enumerate}
\item Under the above assumptions, if $U$ is considered an 
operator of $L^2 (\Omega)\rightarrow L^2 (\Omega)$,then
$U^*U=I$; further if $T$ is injective, then $U$ is the unitary 
operator of $L^2 (\Omega)\rightarrow L^2 (\Omega)$

In fact, since $T$ is measure preserving, $\forall X, Y\in L^2(\Omega)$
\begin{equation*}
\int_ {\Omega} (X\circ T) (Y\circ T) dP=\int_ {\Omega} XY dP
\end{equation*}
Thus $U^*U=I$.

\quad If $T$ is injective, from $P(T(\Omega))=P(\Omega)$, 
for almost every $\omega$, we can define 
$T^{-1}(\omega)$. So for $X\in L^2(\Omega)$, 
we can define $X\circ T^{-1}(\omega)$,
For any $Y\in L^2 (\Omega)$,
\begin{equation*}
\int_ {\Omega} (X\circ T) Y dP=\int_ {\Omega} X (Y\circ T^{-1}) dP
\end{equation*}
From above formula, we see 
$U^*(Y)=Y\circ T^{-1}=U^{-1}(Y)$, that is, $U^*=U^{-1}$, 
that is, $U$ is a unitary operator.

\quad Literature \cite{1}\cite{2}\cite{3} do not 
emphasize the relationship between the 
properties of $T$ (such as injective) and $U$ (such 
as whether it is a unitary transformation). But 
they add some similar conditions, for example
 in Theorem 1 of \cite{1} 
 (i.e. Theorem\ref*{th1} of this paper), 
suppose $E (UU^*\phi|\mathcal{F}_1) = E(\phi|\mathcal{F}_1)$. 
In addition, literature \cite{2} suppose $T$ is sujective in the 
whole text. These conditions are a bit weird at first, 
but as discussed above, they are reasonable. In contrast, 
there are no assumptions about $T$ in literature \cite{7}, 
pointing out directly $U$ it is a unitary transformation.
\item Next, what we need to explain is the relationship 
between the sequence ${U^if}$ in the introduction and 
the stationary, ergodic sequence: If $Y_i=U^if$, 
then the sequence $\{Y_i, i=1\dots\}$ is stationary 
ergodic sequence. Conversely, given a sequence of 
stationary random sequences $\{Y_n, i=1,2, \dots\}$, 
when state space $S$ is good, we can set a probability 
measure $P$ on $S^{N}$ so that the random process 
$\{X_n(\omega)=\omega_ N, n=1,2, \dots\}$ has the same 
distribution as $\{Y_n\}$. For a specific discussion please
see \cite{6} example 6.1.4
\end{enumerate}
\end{rk}
\newpage

\section{Central Limit Theorem of Martingale}
The core method of all literatures reviewed 
is to approximate the sequence ${U^if}$ with 
martingale difference (or backward martingale 
difference). The central limit theorem of 
martingale is introduced below, which is also 
the basis of other results.

\begin{thm}[\cite{4}]
  \textbf{(Martingale's Central Limit Theorem)}\label{lemma1}
  \footnote{Here Thanks Carlanelo 
  Liverani at University of Rome "Tor Vergata"
  Life for helping me find this proof.}
  Let $Y_i$ be a sequence of stationary, ergodic, 
  martingale difference with second-order moments,
  the central limit theorem holds, that is,
  \begin{equation*}
  \frac{1} {\sqrt{n}}\sum_ {i=0}^{n-1} Y_i\rightarrow N (0, \sigma^2)
  \end{equation*}
  Here $\sigma^2=Var(Y_i).$

  \begin{proof}
    First, we define
  \begin{equation*}
    \begin{split}
    &\psi(n,j,t)=exp[\frac{\sigma^2t^2j}{2n}]E\{exp[it\frac{Y_i+\cdots+Y_j}{\sqrt{n}}]\}\\
    &S_n=Y_1+\cdots+Y_n\\
    &\theta(n,j,t)=exp[\frac{\sigma^2t^2j}{2n}]E\{exp[it\frac{S_{j-1}}{\sqrt{n}}]
    [\frac{(\sigma^2-Y_j^2)t^2}{2n}]\}\\
    &\theta_k(n,j,t)=exp[\frac{\sigma^2t^2kr}{2n}]E\{exp[it\frac{S_{kr}}{\sqrt{n}}][\frac{(\sigma^2-Y_j^2)t^2}{2n}]\},\quad kr+1\le j\le k(r+1),
\end{split}
  \end{equation*}
  Here $k$ is an integer.

  To prove the theorem, according to the continuity 
  theorem of characteristic function, it suffices
   to show that for fixed t when $n \rightarrow \infty$,
    we have 
\begin{equation*}
\psi(n,n,t)-1\rightarrow 0
\end{equation*}

Also notes that
\begin{equation*}
\begin{aligned}
|\psi(n,n,t&)-1|=|\sum_ {j=1}^{n}
[\psi(n,j,t)-\psi(n,j-1,t)]|\\&
\le |\sum_ {j=1}^{n}[\psi(n,j,t)-\psi(n,j-1,t)]-
\theta(n,j,t)|+|\sum_ {j=1}^{n}\theta(n,j,t)-
\theta_ k(n,j,t)|\\&+|\sum_ {j=1}^{n}\theta_ k(n,j,t)|
\end{aligned}
\end{equation*}

Next let's estimate the last three 
terms in turn, first of all, for $|t|<T$, we have
\begin{equation}\label{psi-sta}
\begin{aligned}
|[\psi&(n,j,t)-\psi(n,j-1,t)]-\theta(n,j,t)|\\ &=exp[\frac{\sigma^2t^2j}{2n}] 
E\{exp[it\frac{S_{j-1}}{\sqrt{n}}]\{[exp(it\frac{Y_j}{\sqrt{n}})-1-it\frac{Y_j}{\sqrt{n}}+\frac{Y_j^2t^2}{2n}]-
[exp(-\frac{\sigma^2t^2}{2n})-1+\frac{t^2\sigma^2}{2n}]\}\}\\&
\le C(T)E\{|[exp(it\frac{Y_j}{\sqrt{n}})-1-it\frac{Y_j}{\sqrt{n}}+\frac{Y_j^2t^2}{2n}]|\}+C(T)|exp(-\frac{\sigma^2t^2}{2n})-1+\frac{t^2\sigma^2}{2n}|
\end{aligned}
\end{equation}

The first equation uses the properties of martingales:
\begin{equation*}
E\{exp[it\frac{S_{j-1}}{\sqrt{n}}]\xi_ j\}=E(E(exp[it\frac{S_{j-1}}{\sqrt{n}}]\xi_j|S_{j-1}))
=E(exp[it\frac{S_{j-1}}{\sqrt{n}}]E(\xi_j|S_{j-1}))=0
\end{equation*}
$C(T)$ is a constant that only relates to $T$. 
From the formula (\ref{psi-sta})
\begin{equation*}
\sup_{|t|\le T} \sum_ {j=1}^{n}|[\psi(n,j,t)-\psi(n,j-1,t)]-\theta(n,j,t)|=no(\frac{1}{n})
\rightarrow 0
\end{equation*}

Secondly we estimate 
$\sum_ {j=1}^{n}|\theta(n,j,t)-\theta_k (n, j, t)|$, 
select an integer $k$, large and fixed. 
Divide $[1, n]$ into blocks of length k 
in turn. There may be incomplete blocks 
at the end.
\begin{equation}\label{stak-sta}
\begin{aligned}
\sum_ {j=1}^{n}|\theta_ k(n,j,t)-\theta(n,j,t)|&\le
n\sup_ {1\le j\le n} |\theta_ k(n,j,t)-\theta(n,j,t)|\\ & \le C(T)
\sup_ {1\le j\le k} E\{|exp[\frac{\sigma^2t^2j}{2n}]exp[it\frac{S_{j-1}}{\sqrt{n}}]-1||\sigma^2-Y_ j^2|\}
\end{aligned}
\end{equation}
According to the control convergence theorem, for 
each $k$, when $n  \rightarrow \infty$, the above 
equation tends to 0. 

At last we estimate $\sum_ {j=1}^{n}
\theta_k(n, j, t) $, by the stationarity of 
$Y_i$  and $r \le \frac{n}{k}$, we know 
\begin{equation*}
\sum_ {j=kr+1}^{k(r+1)} |\theta_ k(n,j,t)|\le \frac{C(T)}{n} E\{|\sum_{j=kr+1}^{k(r+1)} (\sigma^2-Y_j^2)|\}=C(T)\frac{k}{n}\delta(k)
\end{equation*}

According to Birkhoff's Ergodic Theorem, 
when $k \rightarrow \infty$, $\delta(k) 
\rightarrow 0$,
Since the upper estimate holds for all r, and there 
are at most $ \frac{n}{k} $ blocks, thus we get
\begin{equation*}
\sum_ {j=1}^{n}|\theta_ k(n,j,t)|\le C(T)\delta(k)
\end{equation*}

combining with the formula (\ref{stak-sta}), we have 
\begin{equation*}
\begin{aligned}
|\sum_ {j=1}^{n}\theta(n,j,t)|&\le |\sum_ {j=1}^{n}\theta_ k(n,j,t)|+\sum_ {j=1}^{n}|\theta_ k(n,j,t)-\theta(n,j,t)|
\\& \le C(T)\delta(k)+C(T)\sup_ {1\le j\le k} E\{|exp[\frac{\sigma^2t^2j}{2n}]exp[it\frac{S_{j-1}}{\sqrt{n}}]-1||\sigma^2-Y_ j^2|\}
\end{aligned}
\end{equation*}
Let $n \rightarrow \infty $, then let $k
 \rightarrow  \infty$, we get $ lim_{n\rightarrow \infty}|\sum_ {j=1}^{n}\theta(n,j,t)|\rightarrow 0$

Combining with formula (\ref{psi-sta}), we get
\begin{equation*}
|\sum_ {j=1}^{n}\psi(n,j,t)-\psi(n,j-1,t)|\le |[\psi(n,j,t)-\psi(n,j-1,t)]-\theta(n,j,t)|
+|\sum_ {j=1}^{n}\theta(n,j,t)|\rightarrow 0 
\end{equation*}
\end{proof}
\end{thm}

\section{Backward martingale difference approximation}

Using the martingale central limit theorem 
in the previous section, 
the first main result is proved in this section.
The result is first proved in \cite{1}

\begin{thm}[\cite{1}]\label{th1}
Consider a probability space $(\Omega, \mathcal{F},P)$
and a measure preserving mapping $ T:\Omega \rightarrow 
\Omega$. we assume that if $A\in \mathcal{F}$, 
then $T(A)\in \mathcal{F}$, and the measure preserving
system $(P, T)$ is ergodic. we denote $\mathcal{F}_ 0$
 a sub $\sigma$ algebra of $\mathcal{F}$ and define
$\mathcal{F}_ i=T^{-i}\mathcal{F}_ 0, i\in\mathbb {Z}$, 
and $ \dots  \supset  \mathcal {F}_ {-1}\supset
\mathcal{F}_ 0\supset\mathcal{F}_ 1\supset\dots$ and 
for $ \phi \in L^{\infty}(\Omega) $, we have 
\begin{equation}\label{con}
E(UU^*\phi|\mathcal{F}_1)=E(\phi|\mathcal{F}_1)
\end{equation}

Then for each $f \in L^{\infty}(\Omega)$ which satisfies:\\
(1)$E(f)=0, E(f|\mathcal{F}_0)=f$,\\
(2)$\sum_ {n=0}^{\infty}|E(fU^n f)|< \infty$,\\
(3)$\sum_ {n=0} ^ {\infty} E (U^{*n} f |\mathcal{F}_0)$ converges almost everywhere,
\\Then the central limit theorem is valid, that is, 
\begin{equation*}
\frac{1}{\sqrt{n}} \sum_ {i=0}^{n}U^i f\rightarrow N(0,\sigma^2)
\end{equation*}
Here $\sigma$ satisfies $\sigma ^2  \le - 
E(f^2)+2\sum_ {n=0}^{\infty}|E(fU^n f)|$.

Further, $\sigma=0 $ if and only if there exists
a  $\mathcal{F}_0$ measurable function $g$, such that
\begin{equation*}
Uf=Ug-g
\end{equation*}

Finally, if condition (2) converges in the sense of $L^1$
, then $\sigma ^ 2=- E(f^2)+2\sum_{n=0}^{\infty}E(fU^nf).$
\begin{rk}
Due to the remark 1 of the section 1, $U ^ * U=I$
is always valid. Note that this theorem does not 
assume that $T$ is injective, $U$ is not 
necessarily a unitary transformation, 
so the condition ( \ref*{con}) is 
to ensure that $UU^*=I $ to some extent
\end{rk}

Before formally proving the theorem, we first prove three lemmas
\begin{lem}\label{lem2}
Set a family of $\sigma$ fields $\mathcal 
{F}_ 0\supset\mathcal{F}_ 1  \supset  \dots $
 and a family of random variables $Y_ i, i=1,2 \dotsc $ 
 satisfies that $Y_i ,i \ge 1$ is $F_ {i-1}$measurable,
and
\begin{equation*}
\begin{aligned}
E(Y_i|\mathcal{F}_i)=0,i\ge 0
\end{aligned}
\end{equation*}
Assume also $\sum_ {i=0}^{\infty} Y_i $ converges 
everywhere, then $Y_i$ is backward martingale difference

\begin{proof}
Denote $X_ i=Y_ {i+1}+Y_ {i+2}+\dotsc, i=0,1, \dotsc$
, then $X_i $ is $\mathcal{F}_i $ measurable by the 
 condition
 with
\begin{equation*}
\begin{aligned}
E(X_i|\mathcal{F}_{i+1})=&E(Y_{i+1}+Y_{i+2}+\dotsc|\mathcal{F}_{i+1})\\
=&Y_ {i+2}+Y_ {i+3}\dotsc\\
=&X_ {i+1}
\end{aligned}
\end{equation*}
So we know that $X_i$ is backward martingale difference.
\end{proof}
\end{lem}

\begin{lem}\label{lem3}
Suppose a sequence of random variable $X_n$ converges 
to $Z$, $Y_n$ converges to zero in distribution, 
then $X_ n+Y_ n $ converges to $Z$ in distribution.
\begin{proof}
For any continuous point of the distribution 
function of $Z$, $d$ for example, we have 
\begin{equation*}
P(X_n+Y_n\le d)\le P(X_n\le d+\epsilon)+P(|Y_n|\ge \epsilon)
\end{equation*}
let $n \rightarrow \infty $, and then let $ \epsilon  
\rightarrow 0 $, we get 
$\limsup P (X_n+Y_n  \le d)  \le P (Z  \le d) $\\
For the same reason, from
\begin{equation*}
P(X_n+Y_n\le d)\ge P(X_n\le d-\epsilon)-P(|Y_n|\ge \epsilon)
\end{equation*}
we get 
$ \liminf P (X_n+Y_n  \le d)  \ge P (Z  \le d)$
. That is, $ \lim P (X_n+Y_n \le d)=P(Z \le d)$
. Conclusion follows.
\end{proof}
\end{lem}

\begin{lem}\label{cla1}
$E(U^n\phi|\mathcal{F}_n)=U^nE(\phi|\mathcal{F}_0).$
\end{lem}
This can be obtained from the transformation formula 
of conditional mathematical expectation integral, 
see Section 2.4.7 of the Elements of Probability \cite{5} for details

\begin{proof} [\textbf{proof of theorem \ref*{th1}.}]

  We have three steps to prove the first part
of the theorem \ref*{th1}.\\

\textbf{step1} We want to decompose 
$U ^ nf $ as follows:
\begin{equation}\label{key}
U^n f=Y_ n+U^n g-U^{n-1}g,n\ge 1
\end{equation}
Here, $g$ is almost everywhere finite
 and $\mathcal{F}_ 0$ measurable. $Y_i$ 
 is the backward martingale difference,
  in particular, $Y_i$ satisfies
$E(Y_i|\mathcal{F}_i)=0,i\ge 0$.

Taking $n=1$ in (\ref*{key}), we get 
$Uf=Y_1+Ug-g$,
taking conditional expectation of
 $\mathcal{F}_1$ on both side, 
combining with lemma \ref*{cla1}, 
we obtain $UE(f|\mathcal {F}_0)=
UE(g |\mathcal {F}_0) - E (g | \mathcal {F}_1)$
, notice that $f, g$ is $\mathcal {F}_ 0$ measurable, 
and then act on $U^*$, we get 
\begin{equation}
f=g-U^*E(g|\mathcal{F}_1)=g-U^*E(UU^*g|\mathcal{F}_1)=g-E(U^*g|\mathcal{F}_0)
\end{equation}
Denote $T_ 0: \phi  \mapsto 
E(U ^ *\phi |\mathcal {F}_0) $, 
based on the knowledge of 
functional analysis, the equation
\begin{equation*}
f=(I-T_0)g
\end{equation*}
has a unique solution
 $g=\sum_ {n=0}^{\infty}T_ 0^n f$

\begin{claimm}
$T_ 0^n f=E(U^{*n}|\mathcal{F}_0)$
\end{claimm}
In fact, it suffices to notice that
\begin{equation*}
\begin{aligned}
T_ 0(E(U^{*n}f|\mathcal{F}_0))=&E(U^{*}E(U^{*n}f|\mathcal{F}_0)|\mathcal{F}_ 0)=U^{*}UE(U^{*}E(U^{*n}f|\mathcal{F}_0)|\mathcal{F}_ 0)\\
=&U^{*}E(E(U^{*n}f|\mathcal{F}_0)|\mathcal{F}_ 1)=U^{*}E(U^{*n}f|\mathcal{F}_1)
=U^{*}E(UU^{*(n+1)}f|\mathcal{F}_ 1)\\=& E(U^{*(n+1)}f|\mathcal{F}_ 0)
\end{aligned}
\end{equation*}

So we get
\begin{equation}
g=\sum_ {n=0}^{\infty} E(U^{*n}f|\mathcal{F}_0)
\end{equation}

From the assumptions, we can see 
that $g$ is well defined. Next, we 
will bring the above formula into (\ref{key})
and prove that the decomposition of 
(\ref{key}) is reasonable. 

First, obviously, $g$ is $\mathcal {F}_0$ measurable
 and after calculation, we get 
\begin{equation*}
\begin{aligned}
Y_ i=\sum_ {i=0}^{\infty}E(U^{*n}f|\mathcal{F}_{i-1})-\sum_ {i=0}^{\infty}
E(U^{*n}f|\mathcal{F}_{i}),\quad i\ge 1
\end{aligned}
\end{equation*}
It is easy to know $Y_ i$ is 
$\mathcal {F}_ {i-1} $ measurable, and$
E (Y_i |  \mathcal {F}_i)=0, i \ge 1 $, and 
$\sum_ {i=0}^{\infty} Y_i $ converges almost 
everywhere to $\sum _ {i=0}^{\infty}E(U^{*n}f|\mathcal{F}_0)$
, according to lemma  \ref*{lem2}, $Y_i $ 
is backward martingale difference. 

\begin{claimm}\label{cla3}
$Y_ i$ is a stationary sequence
\end{claimm}
In fact, from the formula (\ref*{key}),
 $Y_ i=U^{i-1}Y_ 1=Y_1 \circ T ^ {i-1}$, 
 where $T$ is the measure preserving 
 transformation. For the positive integer
 $n,m$, and the measurable 
 set $A_ i\in\mathcal{F},i=1\dots n$
we have
\begin{equation*}
\begin{aligned}
P(Y_{m+1}\in A_1,Y_{m+2}\in A_2,\dots,Y_{m+n}\in A_n)&=P(Y_1\in T^{m+1}A_1\cap\dots\cap T^{m+n}A_n)\\ &=P(Y_1\in T^{1}A_1\cap\dots\cap T^{n}A_n)\\&=
P(Y_{1}\in A_1,Y_{2}\in A_2,\dots,Y_{n}\in A_n)
\end{aligned}
\end{equation*}
thus  $Y_i$ is a stationary sequence\\

\textbf{step2}
From the above decomposition,
\begin{equation*}
\frac{1}{\sqrt{n}} \sum_ {i=0}^{n}U^i f=\frac{1}{\sqrt{n}}
\sum_ {i=1}^{n} Y_ i+\frac{1}{\sqrt{n}}(U^n g-g+f)
\end{equation*}
Then $\frac {1}{\sqrt {n}} (U ^ n g-g+f)$ 
converges to 0 in probability. In fact,
for any $\epsilon>0$, and $T$ is a measure
 guaranteed mapping, we have
\begin{equation*}
\begin{aligned}
P(|\frac{1}{\sqrt{n}}(U^n g-g+f)|>\epsilon)=&P(|(U^n g-g+f)|>\epsilon\sqrt{n})\\\le&
P(|(U^n g|>\frac{\epsilon}{2}\sqrt{n})+P(|(f-g|>\frac{\epsilon}{2}\sqrt{n})
\\=&P(|(g|>\frac{\epsilon}{2}\sqrt{n})+P(|(f-g|>\frac{\epsilon}{2}\sqrt{n})
\end{aligned}
\end{equation*}
By $f\in L^{\infty} (\Omega)$, $g$ 
are almost bounded everywhere. It can 
be seen that the above equation is
tend to 0 when $n \rightarrow \infty$.

From lemma  \ref*{lem3}, we know if we 
want to prove 
$\frac {1}{\sqrt {n}} \sum_ {i=0}^{n}U^if$ 
converges in distribution
, it suffices to prove $Y_i$ is a square integrable, 
stationary, ergodic backward martingale difference. 

 That $Y_i$ is stationary is obtained from the claim  
 \ref*{cla3}. Note that $Y_i=U^{i-1}Y_ 1=Y_ 1 \circ
  T ^ {i-1}$, ergodicity is guaranteed by the 
  measure preserving system, and it only needs 
  to prove $Y_i$ is square integrable\\
\textbf{step3}
we will prove $Y_i$ is square integrable in this step. 

In fact, the method is similar to the previous one, 
that is, using martingale difference to approximate 
$Y_i$, we want to find $Y_i(\lambda) $, $ \lambda>1 $
, making $Y_ {i}(\lambda) $ is $\mathcal {F}_{i-1}$ 
measurable, and
\begin{equation*}
\quad  E(Y_{i}(\lambda)|\mathcal{F}_ {i})=0
\end{equation*}
as well as
\begin{equation*}
U^i f=Y_ i(\lambda)+U^ig(\lambda)-\lambda^{-1}U^{i-1}g(\lambda)
\end{equation*}
The same discussion as before shows that
 $g(\lambda)=\sum_ {n=0}^{\infty} 
 \lambda^{-n}E(U^{*n}f|\mathcal{F}_0)$
, note $E (U ^ {* n} f | \mathcal {F}_0)  \le |
 | f||_ {L ^ { \infty}} $ and $\lambda>1$
. So $g (\lambda) \in L_ {\infty}\subset L_ {2} $, and
 $ \lim_ { \lambda \rightarrow 1} g (\lambda)=g (1)=g$, so
$\lim_ {\lambda\rightarrow 1} Y_i(\lambda)=Y_i$.

Due to the stationarity and that U is a
unitary transformation, we have 
\begin{equation*}
\begin{aligned}
E(Uf(Uf-Ug(\lambda)+\lambda^{-1}g(\lambda)))&=E(E(Uf(Uf-Ug(\lambda)+\lambda^{-1}g(\lambda))|\mathcal{F}_ 1))\\ &=
E(UfE((Uf-Ug(\lambda)+\lambda^{-1}g(\lambda))|\mathcal{F}_ 1))\\ &=E(UfE((Y_1|\mathcal{F}_1))=0
\end{aligned}
\end{equation*}
therefore
\begin{equation}\label{long}
\begin{aligned}
E(Y_i(\lambda)^2)=&E(Y_1(\lambda)^2)=E([Uf-Ug(\lambda)+\lambda^{-1}g(\lambda)]^2)\\
=&-E((Uf)^2)+E([Ug(\lambda)-\lambda^{-1}g(\lambda)])\\
=&-E(f^2)+E(Ug(\lambda)[Ug(\lambda)-\lambda^{-1}g(\lambda)])-\lambda^{-1}E(g(\lambda Ug(\lambda)))
+\lambda^{-2}E(Ug(\lambda)^2)\\
=&-E(f^2)+2E(Ug(\lambda)[Ug(\lambda)-\lambda^{-1}g(\lambda)])-(1-\lambda^{-2}E(g(\lambda)^2))\\
=&-E(f^2)+2E(Ug(\lambda)f)-(1-\lambda^2)E(g(\lambda)^2)\\
=&-E(f^2)+2E(g(\lambda)f)-(1-\lambda^2)E(g(\lambda)^2)\\
\le&-E(f^2)+2\sum_ {n=0}^{\infty}\lambda^{-n}E(fU^{*n}f)\\
=&-E(f^2)+2\sum_ {n=0}^{\infty}\lambda^{-n}E(fU^nf)\\
\le&-E(f^2)+2\sum_ {n=0}^{\infty}|E(fU^nf)|
\end{aligned}
\end{equation}
Finally, by the fatous lemma,
\begin{equation*}
E(Y_1^2)=E(\liminf_{\lambda\rightarrow 1} Y_1(\lambda)^2)\le 
\liminf_ {\lambda\rightarrow 1} E(Y_1(\lambda)^2)\le -E(f^2)+2\sum_ {n=0}^{\infty}|E(fU^nf)|
\end{equation*}

Thus $Y_1$ is square integrable. 
So far, the first part of the theorem has been proved\\

For the second part of the theorem, we let $n=1$ 
in the formula (\ref*{key}) and get
\begin{equation*}
\sigma^2=E(Y_1^2)=E[(Uf-Ug+g)^2]
\end{equation*}
Therefore, $\sigma=0  \Longleftrightarrow  
\exists  \mathcal {F}_ 0$
measurable function $g$, such that $Uf=Ug-g$

For the third part of the theorem, note that
\begin{equation}
\begin{aligned}\label{long4}
|E(Y_1^2)-&\{-E(f^2)+2\sum_{n=0}^{\infty}E(fU^nf)\}|\le\\ & |E(Y_1^2)-E(Y_1(\lambda)^2)|+
|E(Y_1(\lambda)^2)-\{-E(f^2)+2\sum_{n=0}^{\infty}E(fU^nf)\}|
\end{aligned}
\end{equation}
We prove that the above two formulas tend to $0$ when 
$\lambda \rightarrow 1 $, so that the conclusion follows. 
For the convenience of proof, 
we first consider the second term, combining with the
 formula \ref*{long}, we have 
\begin{equation}\label{long2}
\begin{aligned}
|E(Y_1(\lambda)^2)-&\{-E(f^2)+2\sum_{n=0}^{\infty}E(fU^nf)\}|\\ &\le 
2\sum_ {n=0}^{\infty} (1-\lambda^{-n})|E(fU^nf)|+(1-\lambda^{-2})E(g(\lambda)^2)\\ &\le
2(1-\lambda^M)\sum_ {n=0}^{\infty} |E(fU^nf)|+ 2\sum_ {n=M}^{\infty} |E(fU^nf)|+(1-\lambda^{-2})E(g(\lambda^2))
\end{aligned}
\end{equation}
Here, $M$ is a large positive integer and fixed. 
To estimate $(1 - \lambda ^ {- 2}) E (g (\lambda ^ 2))$
, for $\lambda>1,\mu>1$,
\begin{equation}\label{long3}
\begin{aligned}
E(g(\lambda)g(\mu))=&\sum_ {n=0,m=0}^{\infty}\lambda^{-n}\mu^{-m}E(U^{*n}fE(U^{*m}f|\mathcal{F}_0))\\ &
\le \sum_ {n=0}^{\infty}\lambda^{-n}\sum_ {m=0}^{M-1} \Vert f\Vert_ {\infty}E(|E(U^{*n}f|\mathcal{F}_0)|)
+\sum_ {n=0}^{\infty}\lambda^{-n}\sum_ {m=M}^{\infty}\Vert f\Vert_ {\infty} E(|E(U^{*m}f|\mathcal{F}_0)|)\\ &
\le M\Vert f\Vert_ {\infty}\sum_ {n=0}^{\infty}E(|E(U^{*n}f|\mathcal{F}_0)|)+\frac{\Vert f\Vert_{\infty}}{1-\lambda^{-1}}\sum_ {m=M}^{\infty} E(|E(U^{*m}f|\mathcal{F}_0)|)
\end{aligned}
\end{equation}
Combining the formula (\ref*{long2}) and 
formula (\ref*{long3}), we can see that 
the second term of (\ref*{long4})
\begin{equation*}
\lim_ {\lambda\rightarrow 1}|E(Y_1(\lambda)^2)-\{-E(f^2)+2\sum_{n=0}^{\infty}E(fU^nf)\}|=0
\end{equation*}
For the first term (\ref*{long4}), notice that 
$\lambda \ge \mu>1 $,
\begin{equation*}
\begin{aligned}
E([Y_1(\lambda)-Y_1(\mu)]^2)=&E([\lambda^{-1}g(\lambda)-\mu^{-1}g(\mu)][Y_1(\lambda)-Y_1(\mu)])\\
=& E([\lambda^{-1}g(\lambda)-\mu^{-1}g(\mu)]^2)+E([g(\lambda)-g(\mu)]^2)\\
\le & (1-\lambda^{-1}\mu^{-1})E(g(\lambda)g(\mu))
\end{aligned}
\end{equation*}
This means that in the sense of $L^2 $,
 $lim_ {\lambda\rightarrow 1}Y_ 1(\lambda)=Y_ 1 $, so
\begin{equation*}
\lim_ {\lambda\rightarrow 1} E(Y_1^2)-E(Y_1(\lambda)^2)=0
\end{equation*}
To sum up, the third part of the theorem is valid
\end{proof}
\end{thm}

For the case that $T$ is a injective, we have the following theorem:
\begin{thm}[\cite{1}]
Assume that $T$ is a invertible  measure preserving 
mapping, $\mathcal {F}_ i\subset \mathcal{F}_ {i-1}$,
 then for $f \in L ^ {\infty} (X), E(f)=0 $. 
 satisfies:\\
(1)$\sum_ {n=0}^{\infty}|E(fU^n f)|< \infty$,\\
(2)$\sum_ {n=0} ^ {\infty} E (U ^ {- n} f | \mathcal {F}_0) $ 
convergence in $L ^ 1$, \\
(3)$\exists\alpha>1:\sup_ {k\in N}k^\alpha 
E(|E(f|\mathcal{F}_ {-k})-f|)<\infty,$

Then the central theorem holds, that is
\begin{equation*}
\frac{1}{\sqrt{n}} \sum_ {i=0}^{n}U^i f\rightarrow N(0,\sigma^2)
\end{equation*}
Here $\sigma ^ 2=- E (f ^ 2)+2 \sum_ {n=0}^{\infty}E(fU^n f).$
\end{thm}
\begin{rk}
The proof idea is similar to the previous theorem. 
The following is just the idea of proof. See \cite{1} for the details
\end{rk}
\begin{proof}
Note that we do not assume that 
$f$ is $\mathcal {F}_ 0$ measurable, so the main idea is to
 use $E (f | \mathcal {F} _ {- k})$ to approximate $f$, 
 that is, using the method in the previous theorem to prove
\begin{equation*}
S_ n^k=\frac{1}{\sqrt{n}}\sum_ {i=0}^{n-1}U^iE(f|\mathcal{F}_{-k})
\end{equation*}
converge to $N (0,  \sigma_k ^ 2) $ in distribution,
 and then let $k  \rightarrow  \infty $, proving that the 
 left side tends to $S_n$, right side tend to 
 $\frac{1}{\sqrt{n}}\sum_ {i=0} ^ {n-1} U ^ if $. And result 
 follows. 
\end{proof}

\newpage
\section{Martingale difference approximation}
This section will be discussed under the condition
 (\ref*{jiashe}) (2) in the section 1, and 
 it is assumed that $T$ is a injective
  (thus, $U$ is a unitary transformation).
  Let us make some notational conventions first. 
  we denote $H_i=L ^ 2(\mathcal {F}_i)$ and
   $S_k=H_ {k+1}\ominus H_k $. we use $Q$ to 
   represent the span of the elements in the form of 
   the $H_ k\ominus H_j $ in $L ^ 2 (\Omega) $, $k$ and $j$ are two integers.
$ P_ {S_k} $ means the projection operator from $L ^ 2 (\Omega)  \rightarrow
 S_k$.

\begin{lem}[\cite{3}]\label{gor}
Set $f \in L ^ 2 (\Omega) $ as well as 
\begin{equation}
\inf_ {g\in Q} \limsup_ {n\rightarrow \infty}^{——} n^{-1}E[\sum_{k=0}^{n-1} U^k(f-g)]^2=0
\end{equation}
then
\begin{equation*}
\lim_ {n\rightarrow\infty} n^{-1} E(\sum_{k=0}^{n-1} U^kf)^2=\sigma^2, 0\le \sigma^2<\infty
,\frac{1}{\sqrt{n}} \sum_ {i=0}^{n}U^i f\rightarrow N(0,\sigma^2).
\end{equation*}
\end{lem}
This lemma is important. The following two 
theorems will be obtained by this lemma. 
 But since the original literature  \cite{3} is 
 written in Russian, which bring many difficulties
  for me to read and understand the proof. The 
  following is my proof by referring to an 
  English document \cite{7} which 
  annotate \cite{3} and some formulas of \cite{3} itself. 

  Before formally proving the theorem, we first prove an assertion
\begin{claimm}\label{cla2}
$H_ k=U^{k}H_ 0.$
\end{claimm}
In fact, according to \ref*{cla1}, 
$\forall f  \in H_ 0 $, $E (U ^ kf | \mathcal {F}_k)
=U ^ kE (f | \mathcal {F}_0)=U ^ kf $, so 
$U ^ kf $ is $ \mathcal {F}_ K $ is measurable.
In addition, $U$ is unitary transformation,
$E(U^kf)^2=E (f ^ 2)<\infty $, which means $U^k f\in H_k$.

On the other hand, $ \forall g \in H ^ k$, 
we let $f=U ^{-n}g $. As above, we have 
$E (f | \mathcal {F}_0)=U ^ {- n} E (U ^ kf | \mathcal {F}_k)=
U ^ {- n} g=f$, so $f $ is $\mathcal {F}_ 0 $ measurable. 
noting $E (f^2)=E(g^2)< \infty $, it means $f\in H_0 $, 
the conclusion is valid.

\begin {proof} [proof of lemma  \ref*{gor}.]
By condition, $p \in \mathbb {Z}^+, 
\epsilon_ p>0, 
\lim_ {p\rightarrow\infty}\epsilon_p=0
 $, there exists$g_p \in Q $ such that 

 \begin{equation*}
\limsup_{n\rightarrow \infty} n^{-1}E[\sum_{k=0}^{n-1} U^k(f-g_p)]^2<\epsilon_ p
\end{equation*}

\textbf{step1} we want to prove the existence of
 random variable $h_ p\in S_{- 1} $ such that
\begin{equation*}
\limsup_{n\rightarrow \infty} n^{-1}E[\sum_{k=0}^{n-1} U^k(f-h_p)]^2<2\epsilon_ p
\end{equation*}
in fact,
\begin{equation*}
\begin{aligned}
f=&g_ p+f-g_ p=\sum_ {l=-\infty}^{\infty} P_ {S_l}g_ p+f-g_ p\\ =&\sum_ {l=-\infty}^{\infty} U^{-(l+1)}P_ {S_l}g_ p+
\sum_ {l=-\infty}^{\infty}\sum_ {m=0}^{-l} U^{m}P_ {S_l}g_ p+U\sum_ {l=-\infty}^{\infty}\sum_ {m=0}^{-l} U^{m}P_ {S_l}g_ p+f-g_ p
\end{aligned}
\end{equation*}
we denote first term of the above formula is $h_p$,
 the second item is $t_p$, then the third item is 
 $Ut_p$, by $g_p\in Q $ we know $t_p$ is well defined, in addition, we have 

\begin{equation*}
\begin{aligned}
\limsup_ {n\rightarrow \infty} n^{-1}E[\sum_{k=0}^{n-1} U^k(f-h_p)]^2=&\lim_ {n\rightarrow \infty}^{——} n^{-1}
E[\sum_{k=0}^{n-1}U^k(t_p-Ut_p+f-g_p)]^2\\ \le& \lim_ {n\rightarrow \infty}^{——} 2n^{-1}[E(t_p-U^nt_p)^2+E(\sum_{k=0}^{n-1}U^k(f-g_p)^2)]<2\epsilon_ p
\end{aligned}
\end{equation*}
To prove the last inequality, it suffices to prove
 $ lim_ {n  \rightarrow \infty} n ^ {- 1} E (t_p-U ^ nt_p) ^ 2=0
  $. From the definition of $Q$ and the claim  
  \ref*{cla2}, we know that when $n$ is sufficiently large, 
  $E (t_pU ^ nt_p)=0 $, as a result, we have 
\begin{equation*}
\begin{aligned}
\lim_ {n\rightarrow \infty} n^{-1}[E(t_p-U^nt_p)^2]=&\lim_ {n\rightarrow \infty} n^{-1}[E(t_p^2)+E(U^nt_p)^2-E(2t_pUt_p)]\\ =&\lim_ {n\rightarrow \infty} n^{-1}[E(t_p^2)+E(t_p)^2]=0
\end{aligned}
\end{equation*}

\textbf{step2} we will prove that 
$ {h_p, p  \in  \mathbb {Z} ^+} $ is 
a convergence sequence of $S_{-1}$ in 
$L^2$.\\

In fact, by $U ^k h_ p\in S_{k-1}$, 
and the orthogonality between $S_k$ and $S_l$
, $k\neq l$,  we get 
\begin{equation*}
\begin{aligned}
E[h_p-h_{p'}]^2&=n^{-1}E[\sum_{k=0}^{n-1}U^k(h_p-h_{p'})]^2=n^{-1}E[\sum_{k=0}^{n-1}U^k(f-h_p)-(f-h_{p'})]^2\\
& \le\limsup_{n\rightarrow\infty}
2n^{-1}E[\sum_{k=0}^{n-1}U^k(f-h_p)]^2+2n^{-1}E[\sum_{k=0}^{n-1}U^k(f-h_{p'})]^2\\
& \le2\epsilon_ p+2\epsilon_ {p'} 
\end{aligned}
\end{equation*}
Thus $ {h_p, p  \in  \mathbb {Z}^+} $ is 
a convergence sequence of $S_ {- 1} $ in 
$L^2$, and $ lim_ {p\rightarrow \infty}h_ p=h_ 0\in
 S_ {-1}$\\

\textbf {step3} we now prove:
\begin{equation*}
\limsup_ {n\rightarrow \infty} n^{-1}E[\sum_{k=0}^{n-1} U^k(f-h_0)]^2=0
\end{equation*}
in fact,
\begin{equation*}
\begin{aligned}
\limsup_ {n\rightarrow \infty} n^{-1}E[\sum_{k=0}^{n-1} U^k(f-h_0)]^2&\le 
\limsup_ {n\rightarrow \infty}  2n^{-1}E[\sum_{k=0}^{n-1} U^k(f-h_p)]^2+n^{-1}E[\sum_{k=0}^{n-1}U^k(h_p-h_0)]^2
\\ & \le \epsilon_ p+E[h_p-h_{p'}]^2
\end{aligned}
\end{equation*}
let $p  \rightarrow \infty$ conclusion follows.\\

\textbf{step4}\\
We denote $r_k=U^k (f-h_0) $, then 
$U ^ kf=U ^ kh_ 0+r_ k$. It is easy to see that 
$ {U ^ kh_0 } $ is the martingale difference sequence 
 \footnote {from $U ^ k h_0 \in S_ {k-1} $,
it can be seen that $E (U ^ kh_0 | \mathcal {F} _ {k-1})=0 $ 
and $E (U ^ kh_0) $ is $ \mathcal {F}_k$ measurable, 
so $U ^ kh_ 0 $ is martingale difference sequence.}, 
 by step 3, we have 
\begin{equation*}
\lim_ {n\rightarrow \infty}^{——} n^{-1}E[\sum_{k=0}^{n-1} r_k]^2=0
\end{equation*}
Therefore, ${r_k}$ converges to 0 in probability, and the 
conclusion is established by the central limit theorem of 
theorem  \ref*{lemma1} martingale and lemma \ref*{lem3}.
\end{proof}
Using this lemma, we prove the following theorems:
\begin{thm}[\cite{2}]
Consider a probability space $(\Omega, \mathcal {F}, P)$
And a measure preserving mapping 
$T: \Omega \rightarrow \Omega $. Here, we also assume that 
$T$ is injective, and if $A \in  \mathcal {F}$,
 then $T (A)  \in  \mathcal {F} $, and further assume that 
$(P, T) $ is ergodic. Let $ \mathcal {F}_ 0 $ be a sub 
$\sigma$ field of $ \mathcal {F} $, defined
$\mathcal{F}_ i=T^{-i}\mathcal{F}_ 0, i \in \mathbb {Z} $, 
and $ \dots  \subset  \mathcal {F}_ {-1}\subset\mathcal{F}_
0\subset\mathcal{F}_ 1  \subset  \dots $, if $f  
\in L ^ 2 ( \Omega) $, denote
\begin{equation*}
x_ r=E(U^rf|\mathcal{F}_0)-E(U^rf|\mathcal{F}_{-1})
\end{equation*}
If $\sum_ {i\in Z} x_ r=Y_ 0\in L^2(\Omega),
E(Y_0^2)=\sigma^2>0$,$
\lim_ {n\rightarrow\infty} n^{-1}ES_n ^ 2= \sigma ^ 2$, \\
then
\begin{equation*}
\frac{1}{\sqrt{n}} \sum_ {i=0}^{n}U^i f\rightarrow N(0,\sigma^2)
\end{equation*}
Here $S_n=\sum_ {i=1}^{n}U^i f.$

\begin{proof}
By condition $Y_0 \in L^2(\Omega) $, 
and $Y_0\in S_ {-1}=H_ 0\ominus H_ {-1}$,
We hope to uses lemma \ref*{gor}.

we Denote $Y_ j=U^j Y_ 0,X_ J=U ^ jf$, and we 
note that 
${Y_j}$ is a stationary ergodic martingale 
difference sequence. \footnote {it can be calculated that 
$U ^ iY_0= \sum_ {r \in \mathbb {Z}}
E [U ^ r f | F_i] - E [U ^ r f | F_ {i-1}] $, 
thus $Y_j $is naturally martingale difference sequence.}, 
If we denote $T_n=\sum_ {i=1}^nY_i $, 
we only need to prove that when 
$n  \rightarrow  \infty$
we have 
\begin{equation*}
\frac{1}{n}E(S_n-T_n)^2\rightarrow 0.
\end{equation*}
Notice 
\begin{equation*}
E(S_n-T_n)^2=ES_ n^2+ET_ n^2-2ES_ nET_ n
\end{equation*}
And by condition and the central limit theorem of 
martingale difference, $n ^ {- 1} ES_ n^2\rightarrow\sigma^2,n^{-1}ET_ n^2=\sigma^2$
, so it is suffices to prove $n ^ {- 1} ES_ nT_ n\rightarrow\sigma^2
$, 
however,
\begin{equation*}
n^{-1}ES_ nT_ n=n^{-1}\sum_ {i=1}^n\sum_ {j=1}^nE(X_iY_j)=\sum_ {j=-(n-1)}^{n-1}(1-|j|n^{-1})E(Y_0X_j)
\end{equation*}
so $n ^ {- 1} ES_ nT_ n\rightarrow\sum_ {j=-\infty}^{\infty}E(Y_0X_j)$
\\On the other hand, $Y_ 0\in H_ 0\ominus H_ {- 1}$
, so $E (Y_0E (X_j | \mathcal {F} _ {- 1}))=0$
, so
\begin{equation*}
E(Y_0X_j)=E(Y_0E(X_j|\mathcal{F}_0))=E(Y_0E(X_j|\mathcal{F}_0)-Y_ 0E(X_j|\mathcal{F}_{-1}))=E(Y_0x_j)
\end{equation*}
So there are
\begin{equation*}
\sum_ {j=-\infty}^{\infty}E(Y_0X_j)=\sum_ {j=-\infty}^{\infty} E(Y_0xj)=E(Y_0^2)=\sigma^2.
\end{equation*}
\end{proof}
\end{thm}
Lemma  \ref*{gor} is a classical result, but its 
condition (\ref*{gor}) is usually difficult to 
verify, so the paper  \cite{2} introduce a 
condition that is stronger than (\ref*{gor}) but 
easier to verify, that is, the following theorem:
\begin{thm}[\cite{2}]
Let $X_j=U ^ jf$, $f$ satisfies $E (f)=0, 
f\in L ^ 2 (\Omega)$ 
and $\mathcal {F}_ 0 $ measurable, 
if it satisfies: \\
(1)$\sum_ {k=1} ^ {\infty} E (X_kE (X_n | \mathcal {F}_0)) $ 
converges for each $n \ge 0 $\\
(2)$\lim_ {n\rightarrow \infty}\sum_ {k=K} ^ {\infty} E (X_kE 
(X_n | \mathcal {F}_0))=0 $ Uniform convergence for $K $\\
Then $ lim_ {n\rightarrow\infty}n^{-1}ES_n ^ 2=\sigma ^ 2, 
\sigma ^ 2< \infty $,if $\sigma ^ 2>0 $, then
\begin{equation*}
\frac{1}{\sqrt{n}} \sum_ {i=0}^{n}U^i f\rightarrow N(0,\sigma^2).
\end{equation*}
\end{thm}
\section {Summary and thanks}

In the previous sections, we listed two different approaches to 
approximate $U^i f$, one is to use backward martingale difference, 
the other is martingale difference. Then we use the martingale 
central limit theorem to obtain the central limit 
theorem about $U ^ if$.

We can also see that the conditions of these theorems are
relatively complex. Whether they have application in 
other fields is what 
I still need to learn and explore.

Finally, I would like to thank Professor Xie Jiansheng for 
his guidance in this semester, and Mr. Carlanelo 
Liverani from the University of Rome "Tor Vergata" for helping 
me find the specific proof of martingale central limit theorem.

\end{document}